\documentclass[12pt]{article}
\usepackage{indentfirst}
\usepackage{amsfonts}
\usepackage[T1]{fontenc}
\usepackage{geometry}
\usepackage{CJK}

\usepackage{amsmath,amssymb,amsthm}
\usepackage{latexsym}
\usepackage{graphicx}
\usepackage{subfigure}
\usepackage{graphicx}
\usepackage[title]{appendix}
\usepackage[square,numbers,sort&compress]{natbib}
\usepackage[colorlinks,
linkcolor=blue,       
anchorcolor=blue,      
urlcolor=blue,         
citecolor=blue,        
]{hyperref}
\usepackage[centerlast]{caption}
\usepackage{times}
\usepackage{subfigure} 
\usepackage{tikz} 
\usepackage{subfigure}
\usepackage{caption}
\usepackage{mathrsfs}
\usepackage{comment}
\usepackage[babel]{microtype}
\usepackage[english]{babel}
\usepackage[capitalise]{cleveref}

\usepackage{enumerate}
\usepackage{appendix}
\usepackage[shortlabels]{enumitem}
\usepackage{graphicx}
\usetikzlibrary{shapes.geometric}
\usepackage{rotating}
\usepackage{booktabs}
\usepackage{xcolor}
\definecolor{revred}{RGB}{190,0,0}

\usetikzlibrary{calc,decorations.pathreplacing,decorations.markings,shapes.geometric}
\tikzstyle{vertice}=[circle,fill,inner sep=1.5pt]
\tikzset{
	on each segment/.style={
		decorate,
		decoration={
			show path construction,
			moveto code={},
			lineto code={
				\path [#1]
				(\tikzinputsegmentfirst) -- (\tikzinputsegmentlast);
			},
			curveto code={
				\path [#1] (\tikzinputsegmentfirst)
				.. controls
				(\tikzinputsegmentsupporta) and (\tikzinputsegmentsupportb)
				..
				(\tikzinputsegmentlast);
			},
			closepath code={
				\path [#1]
				(\tikzinputsegmentfirst) -- (\tikzinputsegmentlast);
			},
		},
	},
	mid arrow/.style={postaction={decorate,decoration={
				markings,
				mark=at position .6 with {\arrow[scale=1]{stealth}}
	}}},    
	lid arrow/.style={postaction={decorate,decoration={
				markings,
				mark=at position .4 with {\arrow[scale=1]{stealth}}
	}}},
	llid arrow/.style={postaction={decorate,decoration={
				markings,
				mark=at position .53 with {\arrow[scale=1]{stealth}}
	}}},
	rid arrow/.style={postaction={decorate,decoration={
				markings,
				mark=at position .7 with {\arrow[scale=1]{stealth}}
	}}},      
}

\newtheorem{theorem}{Theorem}[section]
\newtheorem{proposition}[theorem]{Proposition}

\newtheorem{conjecture}[theorem]{Conjecture}

\newtheorem{lemma}[theorem]{Lemma}
\newtheorem{claim}[theorem]{Claim}

\begin{document}
	\title{Two-block paths in oriented graphs under the Ore-type condition}

	\author{Ke Zhao \quad Jia Zhou  \: Yunshu Gao\thanks{Corresponding author: gysh2004@gmail.com}\\
		School of Mathematics and Statistics, Ningxia University\\
		Yinchuan, China 750021}

	\maketitle

	\bibliographystyle{plain}

	\begin{abstract}

		Stein conjectured that for every positive integer $k$, every oriented graph with minimum semi-degree greater than $k/2$ contains every oriented path of length $k$ [Springer Optim.Appl. 165 (2020), 459--486].   For directed paths this is a classical result of Jackson (1981). The two-block case of the conjecture was established by Chen, Hou and Zhou (2025). In this paper we prove the two-block analogue of Jackson's theorem. 
		For all positive integers $s$ and $r$, let $P(s,r)$ be the oriented path consisting of two directed paths of lengths $s$ and $r$ with a common terminal vertex, and let $Q(s,r)$ be its reverse. We show that for all $K\ge3$ and all $s,r\ge1$ with $s+r\le K$, every oriented graph $D$ with
		$\delta(D)\ge K$ and $d_D^+(u)+d_D^-(v)\ge K-1$ for every ordered pair of distinct vertices
		$(u,v)\notin A(D)$ contains both $P(s,r)$ and $Q(s,r)$. Our result is the two-block analogue of Jackson's theorem, and it extends the two-block result of Chen, Hou and Zhou. The proof iterates a shift of the terminal vertex of a longest path.

	\end{abstract}

	\section{Introduction}

	Throughout this paper, all digraphs considered have no loops or parallel arcs. Let $D=(V,A)$ be a digraph with vertex set $V=V(D)$ and arc set $A=A(D)$. For an arc $(u,v)\in A(D)$, we also say that $u$ \textbf{dominates} $v$. Two different vertices are \textbf{adjacent} if there exists an arc connecting them. For $v\in V(D)$, write $d_D^+(v)$ and $d_D^-(v)$ for the out-degree and in-degree of $v$, and let
	$d_D(v)=d_D^+(v)+d_D^-(v)$ and $ \delta(D)=\min_{v\in V(D)}d_D(v).
	$
	Let $\delta^+(D)=\min\{d_D^+(v)\mid v\in V(D)\}$ and $\delta^-(D)=\min\{d_D^-(v)\mid v\in V(D)\}$  
	be the \textbf{minimum out-degree} and \textbf{minimum in-degree} of $D$. We write $\delta^0(D)=\min\{\delta^+(D),\delta^-(D)\}$ for the \textbf{minimum semi-degree} of $D$.
	An \textbf{oriented graph} is a digraph with no cycles of length two.   For every oriented graph $D$, if  $(u,v)\in A(D)$ then we have $(v,u)\notin A(D)$. A \textbf{block} of an oriented path is a maximal directed
	sub-path of it, thus a two-block path consists of two directed subpaths in opposite directions, which share exactly one common vertex.

	Degree conditions forcing paths and cycles have a long history. Ore~\cite{Ore1960} proved that every
	graph on $n\ge3$ vertices is Hamiltonian if $d(u)+d(v)\ge n$ for every two nonadjacent vertices $u,v$. Woodall~\cite{Woodall1972} carried this Ore-type condition over to digraphs, proving that a digraph $D$ on $n\ge2$ vertices is Hamiltonian if $d_D^+(u)+d_D^-(v)\ge n$ for every ordered pair of distinct vertices $u,v$ with $(u,v)\notin A(D)$. The semi-degree counterpart is due to
	Ghouila-Houri~\cite{GhouilaHouri1960}. For oriented graphs, Keevash, K\"uhn and
	Osthus~\cite{KeevashKuhnOsthus2009} showed that minimum semi-degree at least $(3n-4)/8$ forces Hamiltonicity for all sufficiently large $n$.
	The existence of oriented paths has also been studied extensively. Jackson~\cite{Jackson1981} proved the following theorem.
	\begin{theorem}[\cite{Jackson1981}]\label{thm:jackson}
		Let $n$ be a positive integer. If an oriented graph $D$ satisfies $\delta(D)\ge n$ and $d_D^+(u)+d_D^-(v)\ge n-1$ for every $(u,v)\notin A(D)$, then $D$ contains a directed path of length $n$. 	
	\end{theorem}
	Stein~\cite{Stein2020} proposed the following conjecture.

	\begin{conjecture}[\cite{Stein2020}]\label{Stein}
		For any positive integer \(k\), every oriented graph \(D\) of \(\delta^0(D) > k/2\) contains every oriented path of length \(k\).
	\end{conjecture}
	
	Theorem \ref{thm:jackson} implies that Conjecture \ref{Stein} holds for directed paths. An oriented path is \textbf{antidirected} if no two consecutive arcs are both forward and no two consecutive arcs are both backward. For antidirected paths, Klimo\v{s}ov\'a and Stein~\cite{KlimosovaStein2023} proved that minimum semi-degree at least $(3k-2)/4$ suffices; this bound was improved to $2k/3$ by Chen, Hou and
	Zhou~\cite{ChenHouZhou2024} and further to $5k/8$ by Skokan and
	Tyomkyn~\cite{SkokanTyomkyn2025}. In 2025, Chen, Hou and Zhou~\cite{ChenHouZhou2025} proved  Conjecture \ref{Stein} for two-block paths.

	\begin{theorem}[\cite{ChenHouZhou2025}]\label{Chen}
		For any integer \(k\), every oriented graph \(D\) of \(\delta^0(D) > k/2\) contains every oriented path of length \(k\) with two blocks.
	\end{theorem} 
	
	Theorem \ref{Chen} in turn implies the result of Penev et al.~\cite{PenevEtAl2026}. We mention one more line of research. 
	A \textbf{coloring} of a digraph is a partition of its vertex set into parts that induce no arc,
	and the \textbf{chromatic number} $\chi(D)$ is the least size of such a partition.
	Two-block paths were first studied under chromatic conditions: El Sahili~\cite{ElSahili2004} conjectured that every $k$-chromatic digraph contains every two-block path with $k$ arcs, a conjecture confirmed by Addario-Berry et al.~\cite{AddarioBerryHavetThomasse2007}, see also~\cite{ElSahiliKouider2007}. 
	
	Motivated by these results, we prove the following two-block analogue of Theorem \ref{thm:jackson} under the Ore-type condition. 
	For positive integers $s$ and $r$, let $P(s,r)$ be the oriented path
	$z_0\to z_1\to\cdots\to z_s\leftarrow z_{s+1}\leftarrow\cdots\leftarrow z_{s+r}$. It consists of two internally disjoint directed paths of lengths $s$ and $r$ that share only their common
	terminal vertex $z_s$. Let $Q(s,r)$ be the reverse of $P(s,r)$. Since the two arms of $P(s,r)$ are
	interchangeable, $P(s,r)\cong P(r,s)$ and $Q(s,r)\cong Q(r,s)$. We may therefore
	assume that $s\ge r$ in the proof.
	\begin{theorem}\label{thm:main}
		Let $K\ge3$ and $s, r\ge1$ satisfy $s+r\le K$. Every oriented graph $D$  satisfying $\delta(D)\ge K$ and $d_D^+(u)+d_D^-(v)\ge K-1$ for every $(u,v)\notin A(D)$ contains both $P(s,r)$ and $Q(s,r)$.
	\end{theorem}
	
	The bound $K\ge3$ is needed only when $s=r=1$, since $s,r\ge1$ and $s+r\le K$ force $s=r=1$ when $K=2$. A directed cycle $C$ on $n\ge3$ vertices satisfies $\delta(C)=2$, and $d_C^+(u)+d_C^-(v)=2\ge1$ for every pair $(u,v)\notin A(C)$. Yet no vertex of $C$ has two in-neighbours or two out-neighbours. Consequently, $C$ contains neither $P(1,1)$ nor $Q(1,1)$.

	Theorem~\ref{thm:main} is the two-block analogue of Theorem~\ref{thm:jackson}. Semi-degree bounds form a special case: if $\delta^0(D)\ge k$, then $\delta(D)\ge 2k\ge k$ and $d_D^+(u)+d_D^-(v)\ge 2k>k-1$ for every pair $(u,v)\notin A(D)$. The hypothesis of Theorem~\ref{thm:main} is thus strictly weaker than its semi-degree counterpart; Proposition~\ref{prop:no-semi} imposes no lower bound on $\delta^0(D)$. Consequently, Theorem~\ref{thm:main} generalizes Theorem~\ref{Chen}: for $k\ge3$, if $\delta^0(D)>k/2$, then $\delta^0(D)\ge\lfloor k/2\rfloor+1$, so $\delta(D)\ge k$, $d_D^+(u)+d_D^-(v)\ge k+1>k-1$, and applying Theorem~\ref{thm:main} with $K=k$ yields all $P(s,r)$ with $s+r\le k$.

	\begin{proposition}\label{prop:no-semi}
		For every $K\ge2$ there is an oriented graph $D$ with $\delta(D)\ge K$ and $d_D^+(u)+d_D^-(v)\ge K-1$ for every $(u,v)\notin A(D)$, whose minimum semi-degree is zero.
	\end{proposition}
	
	\begin{proof}
		Let $T$ be a regular tournament on $2K-1$ vertices and let $D$ be obtained from $T$ by adding a
		new vertex $u$ that dominates every vertex of $T$. Every vertex $v\in V(T)$ satisfies
		$d_D^+(v)=K-1$ and $d_D^-(v)=K$, whereas $d_D^+(u)=2K-1$ and $d_D^-(u)=0$, hence
		$d_D(v)=2K-1=d_D(u)$ for every $v\in V(T)$, and $\delta(D)=2K-1\ge K$.
		Consequently, for every pair of vertices $a,b\in V(D)$ with $(a,b)\notin A(D)$, exactly one of the two equalities holds: either $d_D^+(a)+d_D^-(b)=2K-1$, or $d_D^+(a)+d_D^-(b)=K-1$. In both cases, $d_D^+(a)+d_D^-(b)\ge K-1$.
		Since $d_D^-(u)=0$, the minimum semi-degree of $D$ is zero, which completes the construction.
	\end{proof}
	
	The proof uses the theorem of Havet and Thomass\'e~\cite{HavetThomasse2000} on oriented Hamiltonian paths in tournaments.
	
	\begin{theorem}[\cite{HavetThomasse2000}]\label{Havet}
		Let $T$ be a tournament of order $n$ and let $P$ be an oriented path of order $n$. Then $T$
		contains $P$ unless $P$ is antidirected and $T$ is the $3$-cycle, the regular tournament on
		five vertices, or the Paley tournament on seven vertices.
	\end{theorem}

	\noindent \textbf{Organization.}  
	The rest of the paper is organized as follows. Section 2 sets up the notation used throughout the paper. Section 3 gives the proof of Theorem~\ref{thm:main}. Finally, Section 4 concludes the paper with some remarks.

	\section{Notation}
	Undefined notation follows 
	\cite{BangJensenGutin2009}. For integers $a\le b$
	we set $\boldsymbol{[a,b]} =\{a,a+1,\ldots,b\}$. 
	Given a digraph $D$ and a vertex $v \in V(D)$,
	we write $N_D^+(v)$, $N_D^-(v)$ and
	$N_D(v)$  for the  \textbf{out-neighbourhood}, \textbf{in-neighbourhood} and \textbf{neighbourhood} of $v$. 
	For $S \subseteq V(D)$, 
	we set
	$$
	d_D^+(v;S):=|N_D^+(v)\cap S|,\qquad d_D^-(v;S):=|N_D^-(v)\cap S|.
	$$
	The subscript $D$ is omitted whenever no ambiguity arises. 
	We denote by $D[S]$ the sub-digraph of $D$
	induced by $S$, and we call the digraph obtained from $D$ by reversing every
	arc its \textbf{reverse}.
	For $X\subseteq V(D)$, let $D\setminus X$ denote the digraph obtained from $D$ by deleting the vertices in $X$ and all arcs incident to any vertex of $X$.

	A directed path $P=x_1x_2\cdots x_n$ has \textbf {length} $n-1$ and \textbf {order} $n$, with initial and terminal vertices $x_1$ and $x_n$, respectively.  Throughout, indices assigned to vertices of a path always refer to the labelling of that path.

	An \textbf{oriented path} is an orientation of an undirected path. For an oriented path $v_1v_2\cdots v_n$, an arc $(v_i,v_{i+1})$ is called
	\textbf{forward} and an arc $(v_{i+1},v_i)$ is called \textbf{backward}.

	A directed cycle containing every vertex of $D$
	is a \textbf{directed Hamiltonian cycle} of $D$,
	and indices of a directed cycle are always
	taken modulo its length. The digraph $D$ \textbf{contains} $H$ if $H$ is a sub-digraph of $D$.

	A digraph is \textbf{complete} if $(u,v)\in A(D)$ and $(v,u)\in A(D)$ for any two distinct vertices $u,v$. 
	A \textbf{tournament} is an orientation of a complete graph, that is, an oriented graph in which every pair of distinct vertices is adjacent. A tournament is \textbf{regular} if every vertex has the same out-degree.
	We define $P(s,r)$ as the oriented path consisting of $s$ forward arcs followed by $r$ backward arcs, and $Q(s,r)$ as its reverse.

	\section{Proof of Theorem~\ref{thm:main}}
	
	In this section we present the proof of Theorem \ref{thm:main}. Throughout this section, suppose that $D$ satisfies $\delta(D)\ge K$ and
	$d_D^+(u)+d_D^-(v)\ge K-1$ for every $(u,v)\notin A(D)$, where $K\ge3$. First let $s=r=1$. 
	If $d_D^-(v)\le1$ for every $v\in V(D)$ then $|A(D)|\le|V(D)|$, whereas $\delta(D)\ge3$ gives $2|A(D)|\ge3|V(D)|$,
	a contradiction. Hence some vertex has two in-neighbours, giving $P(1,1)$. 
	The same counting with out-degrees instead of in-degrees gives a vertex with two out-neighbours, i.e.\ $Q(1,1)$.
	From now on, we assume that $s+r\ge 3$, $s\ge 2$ and $s\ge r$. We only present the proof of Theorem~\ref{thm:main} when $s+r=K$, since the case $s+r<K$ can be deduced from the case $s+r=K$. 
	We only present the existence of $P(s,r)$,
	while the existence of $Q(s,r)$ follows by reversing all arcs of $D$.

	\subsection{Outline of the proof}
	We first sketch the proof. The proof of Theorem~\ref{thm:main} proceeds by the contradiction: suppose that $D$ contains no $P(s,r)$.
	
	Using Lemma \ref{longest path}, we characterize the distribution of the   neighbourhood of the terminal vertex $x_n$ of $P$, where $P$ is an arbitrary longest path of $D$. Based on this property, the proof proceeds by induction on $r$.
	
	For the base case $r=1$: since $s+r\ge 3$, we have $s\ge 2$. Lemma \ref{longest path} together with Theorem \ref{Havet} implies that $D$ contains $P(s,1)$, which yields a contradiction.
	Now suppose $r\ge 2$. Assume that for every integer $f$ with $1\le f<r$, every oriented graph $D'$ satisfying $\delta(D')\ge s+f$ and
	$d^+_{D'}(u)+d^-_{D'}(v)\ge s+f-1$ for every $(u,v)\notin A(D')$
	contains $P(s,f)$.
	We prove Lemma \ref{lem:shift}, which asserts that for every longest path $P$ in $D$, one may modify the terminal vertex of $P$ to yield a new longest path with the same vertex set $V(P)$. 
	This yields a cyclic longest-path rearrangement algorithm to obtain two disjoint vertex sets $W$ and $U$ such that $W \to U$, $D[W]$ contains the path $P_W:=x_1x_2\cdots x_h$, and $D[U]$ satisfies $\delta(D[U])\ge s+r-h$ and
	$d^+_{D[U]}(u)+d^-_{D[U]}(v)\ge s+r-h-1$ for every $(u,v)\notin A(D[U])$. 
	By the induction hypothesis, $D[U]$ contains $P(s,r-h)$. 
	Attaching the path $P_W$ to this path yields the desired $P(s,r)$, giving a contradiction.

	\subsection{Preliminary lemmas}
	In this section, we prove several lemmas that will be used in the proof. For every longest directed path $P=x_1x_2\cdots x_n$ of $D$, let 
	$R(P)=\{x_{s+1},\ldots,x_{n-r-1}\}$ and $B(P)=\{x_{r+1},\ldots,x_{n-s-1}\}.$
	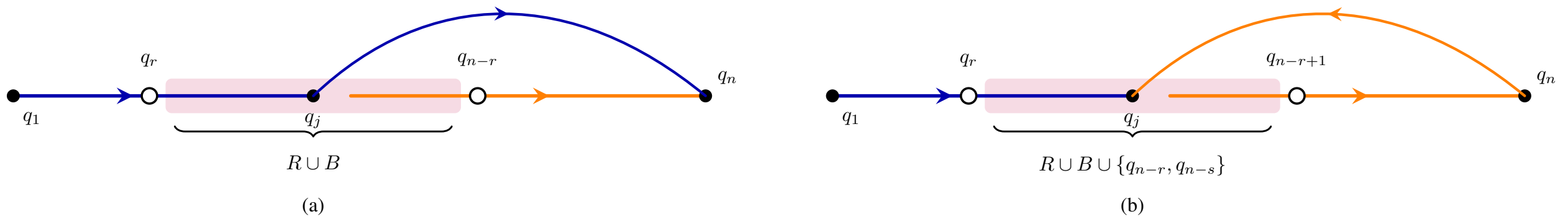
\begin{figure}[htbp]
		\centering
		
		\newcommand{\figTwoBase}[1]{%
			\coordinate (L)  at (0,0);
			\coordinate (T1) at (2.8,0);
			\coordinate (Vi) at (7.6,0);
			\coordinate (R)  at (12.0,0);
			\coordinate (Xj) at (5.2,0);
			\coordinate (G0) at (5.85,0);  
			\fill[purple!14,rounded corners=3pt]
			($(T1)+(-0.16,-0.30)$) rectangle ($(Vi)+(0.16,0.30)$);
			
			\draw[line width=1.7pt,blue!68!black,figlineL] (L) -- (Xj);
			
			\draw[line width=1.7pt,orange,figlineR]        (G0) -- (R);
			
			\fill (L)  circle (3.2pt);
			\fill (R)  circle (3.2pt);
			\fill (Xj) circle (3.2pt);
			
			\node[figlab,anchor=north west] at ($(L)+(0.02,-0.14)$)  {$q_1$};
			\node[figlab,anchor=south west] at ($(R)+(0.06,0.06)$)   {$q_n$};
			\node[figlab,anchor=north]      at ($(Xj)+(0,-0.16)$)    {$q_j$};
			
			\fill[white] (2.36,0) circle (3.8pt);
			\draw[line width=1.1pt] (2.36,0) circle (3.8pt);
			\fill[white] (8.05,0) circle (3.8pt);
			\draw[line width=1.1pt] (8.05,0) circle (3.8pt);
			\node[figlab,anchor=south] at (2.36,0.36) {$q_{r}$};
			\node[figlab,anchor=south] at (8.05,0.36) {$q_{#1}$};
		}
		\newcommand{\figTwoBrace}[1]{%
			\draw[decorate,decoration={brace,amplitude=5pt,mirror,raise=2pt},line width=0.9pt]
			($(T1)+(0,-0.46)$) -- ($(Vi)+(0,-0.46)$);
			\node[figlab,anchor=north] at ($(T1)!0.5!(Vi)+(0,-0.89)$) {#1};
		}
		
		\newcommand{\figTwoArcFwd}{%
			\draw[line width=1.3pt,blue!68!black,figmid]
			(Xj) .. controls ($(Xj)+(2.0,1.9)$) and ($(R)+(-2.3,1.9)$) .. (R);
		}
		
		\newcommand{\figTwoArcBack}{%
			\draw[line width=1.3pt,orange,figmid]
			(R) .. controls ($(R)+(-2.3,1.9)$) and ($(Xj)+(2.0,1.9)$) .. (Xj);
		}
		\resizebox{\textwidth}{!}{%
			\begin{tikzpicture}[line cap=round,line join=round,
				figlab/.style={font=\footnotesize},
				figlineL/.style={postaction={decorate,decoration={markings,
							mark=at position 0.40 with {\arrow[scale=1.30]{stealth}}}}},
				figlineR/.style={postaction={decorate,decoration={markings,
							mark=at position 0.56 with {\arrow[scale=1.30]{stealth}}}}},
				figmid/.style={postaction={decorate,decoration={markings,
							mark=at position 0.5 with {\arrow[scale=1.30]{stealth}}}}}]
				
				\begin{scope}
					\figTwoBase{n-r}
					\figTwoArcFwd
					\figTwoBrace{$R\cup B$}
					\node[figlab] at (5.2,-1.92) {(a)};
				\end{scope}
				
				\begin{scope}[shift={(14.2,0)}]
					\figTwoBase{n-r+1}
					\figTwoArcBack
					\figTwoBrace{$R\cup B\cup\{q_{n-r},q_{n-s}\}$}
					\node[figlab] at (5.2,-1.92) {(b)};
				\end{scope}
		\end{tikzpicture}}
		\caption{Construction of the path $P(s,r)$.}
		\label{fig:construction}
	\end{figure}
	
	\begin{lemma}
		\label{longest path}
		Let $s,r,K$ be positive integers with $s\geq r \geq 1$ and $s+r=K\geq 3$. Suppose that $D$ is an oriented graph with $\delta(D)\ge K$ and
		$d_D^+(u)+d_D^-(v)\ge K-1$ for every $(u,v)\notin A(D)$, but $D$ contains no $P(s,r)$. Suppose $P=x_1x_2\cdots x_n$ is a longest directed path of $D$. Then for every longest directed path $Q=q_1\cdots q_n$  on $V(P)$, the following statements hold.
		\begin{itemize}
			\item[(1)]  $N_D(q_n)\cap(R(Q)\cup B(Q))=\emptyset$. Moreover, $(q_n,q_{n-r})\notin A(D)$ and $(q_n,q_{n-s})\notin A(D)$.
			\item[(2)]
			The terminal vertex $q_n$ has no
			in-neighbour in $D\setminus V(P)$.
			\item[(3)] 
			Put $c=n-s-r-1$, then $ |R(Q)\cup B(Q)|=c.$ Hence $R(Q)=B(Q)$, which forces $s=r$ whenever $c>0$; moreover, if $s\ne r$ then $R(Q)=B(Q)=\emptyset.$
			\item[(4)] 
			$q_n$ is adjacent to every vertex of $V(P)\setminus(B(Q)\cup\{q_n\})$.
		\end{itemize}

	\end{lemma}
	\begin{proof}
		By Theorem~\ref{thm:jackson}, $n\ge K+1=s+r+1$.
		The definition of $R(P)$ and $B(P)$ implies that $|R(P)|=|B(P)|=n-s-r-1$.

		We first prove (1) and only prove the assertion for $R(Q)$, the proof for $B(Q)$ follows by a similar argument. Suppose the assertion is false. Then there exists
		$i\in[s+1,n-r-1]$ such that $q_i$ and $q_n$ are adjacent.  If $(q_i,q_n)\in A(D)$, then
		$$
		q_{i-s+1}\to\cdots\to q_i\to q_n\leftarrow q_{n-1}\leftarrow\cdots\leftarrow q_{n-r}
		$$
		is $P(s,r)$,  a contradiction (see  Figure~\ref{fig:construction} (a)). 
		If $(q_n,q_i)\in A(D)$, then
		$$
		q_{i-s}\to\cdots\to q_i\leftarrow q_n\leftarrow q_{n-1}\leftarrow\cdots\leftarrow q_{n-r+1}
		$$
		is $P(s,r)$, a contradiction (see  Figure~\ref{fig:construction} (b)). Hence no vertex of $R(Q)$ is adjacent to $q_n$.
		
		Finally, suppose that $(q_n, q_{n-r})\in A(D)$. Then
		$$
		q_{n-s-r}\to\cdots\to q_{n-r}\leftarrow q_n\leftarrow q_{n-1}\leftarrow\cdots\leftarrow q_{n-r+1}
		$$
		is $P(s,r)$, a contradiction. The exclusion $(q_n, q_{n-s})\not\in A(D)$ follows by symmetry.
		
		Next, we prove (2). 
		Suppose that $y_1\notin V(P)$, $(y_1,q_n)\in A(D)$. We inductively construct distinct vertices $y_1,\ldots,y_r$ outside $V(P)$ such that
		\begin{equation}
			y_j\to y_{j-1}\to\cdots\to y_1\to q_n \qquad (1\le j\le r). \label{eq:chain}
		\end{equation}
		Assume that $y_1,\ldots,y_j$ have been constructed, $1\le j\le r-1$. 
		Since $Q$ is longest, $q_n$ has no out-neighbour outside $V(P)$, which implies $q_n\not\to y_j$. Together with $d_D^+(q_n)+d_D^-(y_j)\ge s+r-1$,
		$d_D^-(y_j)\ge s+r-1-d_D^+(q_n)$.

		Now we consider the distribution of  the in-neighbours of $y_j$ on $V(P)$.  If there is an index $i$ with $r-j\le i\le n-s-1$ and $(q_i,y_j)\in A(D)$, then
		$$
		q_{i-r+j+1}\to\cdots\to q_i\to y_j\to\cdots\to y_1\to q_n
		$$
		together with $q_{n-s}\to\cdots\to q_n$ forms $P(s,r)$, a contradiction. Thus $y_j$ has no in-neighbour  $q_i$ with $r-j\le i\le n-s-1$.
		Also $(q_{n-1},y_j)\not\in A(D)$, since otherwise
		$$
		q_1\to\cdots\to q_{n-1}\to y_j\to\cdots\to y_1\to q_n
		$$
		would be longer than $P$. Thus the in-neighbours of $y_j$ on $V(P)$ lie in
		\begin{equation}
			A_j=\{q_i\mid i\in [1,r-j-1]\cup[n-s,n-2] \},\qquad |A_j|=s+r-j-2. \label{eq:Aj}
		\end{equation}
		Hence $d_D^-(y_j;A_j)\le s+r-j-2$.
		
		Our next goal is to prove that $y_j$ has at least $j$ in-neighbours outside $V(P)$. 
		Let $q_k$ be an out-neighbour of $q_n$. Then $k>j$, since otherwise
		$$
		y_j\to\cdots\to y_1\to q_n\to q_k\to\cdots\to q_{n-1}
		$$
		would be longer than $P$.  
		Next we consider three cases: $k\le r-1$, $k=r$ and $k\ge r+1$. 
		If $k\le r-1$, then $k-j\in[1,r-j-1]$.
		Further, we obtain $(q_{k-j},y_j)\not\in A(D)$
		since otherwise the path
		$$
		q_1\to\cdots\to q_{k-j}\to y_j\to\cdots\to y_1\to q_n\to q_k\to\cdots\to q_{n-1}
		$$
		would be longer than $P$. 
		When $k=r$, neither $q_{r-j}$ nor $q_{r-1}$ lies in $A_j$, so this case yields no vertex of $A_j$. We therefore consider only $k\le r-1$ and $k\ge r+1$.
		By (1), $k\ge n-s+1$ (see Figure~\ref{fig:construction}). Now we show $(q_{k-1},y_j)\not\in A(D)$, as otherwise 	
		$$
		q_1\to\cdots\to q_{k-1}\to y_j\to\cdots\to y_1\to q_n\to q_k\to\cdots\to q_{n-1}
		$$
		would be longer than $P$. 
		Hence at least $d_D^+(q_n)-1$ vertices of $A_j$ are not in-neighbours of $y_j$. Therefore
		\begin{equation}
			d_D^-(y_j;V(P))\le (s+r-j-2)-(d_D^+(q_n)-1)=s+r-j-d_D^+(q_n)-1, \label{eq:yS}
		\end{equation}
		where for $d_D^+(q_n)=0$ the bound is immediate from~\eqref{eq:Aj}. Combining $d_D^-(y_j)\ge s+r-1-d_D^+(q_n)$ and~\eqref{eq:yS}, we have
		\begin{equation}
			d_D^-(y_j;V\setminus V(P))\ge (s+r-1-d_D^+(q_n))-(s+r-j-d_D^+(q_n)-1)=j. \label{eq:y-out}
		\end{equation}
		
		Finally we construct the path 	
		$y_r\to y_{r-1}\to\cdots\to y_1
		$ by induction. 
		Since $(y_j,y_{j-1})\in A(D)$ and $D$ is oriented, $y_{j-1}$ is not an in-neighbour of $y_j$ (vacuously for $j=1$). Hence among the used vertices $y_1,\ldots,y_{j-1}$, at most $j-2$ are in-neighbours of $y_j$. By~\eqref{eq:y-out}, $y_j$ has at least $j$ in-neighbours outside $V(P)$, one of them is unused, and we take it as $y_{j+1}$. By the induction, there exists the desired path	$y_r\to y_{r-1}\to\cdots\to y_1
		$.  The path
		$
		y_r\to y_{r-1}\to\cdots\to y_1\to q_n
		$
		together with $q_{n-s}\to\cdots\to q_n$ forms $P(s,r)$, a contradiction.
		
		From now on, we prove (3).
		Since $\delta(D)\ge s+r$ and $V(P)\setminus(R(Q)\cup B(Q)\cup\{q_n\})$
		has $n-1-|R(Q)\cup B(Q)|$ vertices,   the statement (1) gives
		$$
		s+r\le d_D(q_n)\le n-1-|R(Q)\cup B(Q)|\le n-1-c=s+r,
		$$
		so equality holds throughout. In particular $d_D(q_n)=s+r$. Since $c\ge0$ and $|R(Q)|=n-s-r-1=c$, we have $c=|R(Q)|\le|R(Q)\cup B(Q)|\le c$, so $|R(Q)|=|B(Q)|=|R(Q)\cup B(Q)|=c$. 
		This yields \(s=r\) whenever \(c>0\), and hence \(c=0\) whenever \(s>r\), which completes the proof of (3).

		For (4), since $|B(Q)|=c$ we have
		$|V(P)\setminus(B(Q)\cup\{q_n\})|=n-c-1=s+r$. Statement (2) and the maximality of $Q$ imply that $q_n$ has no neighbour in $(D\setminus V(P))\cup(B(Q)\cup\{q_n\})$, while
		$d_D(q_n)=s+r=|V(P)\setminus(B(Q)\cup\{q_n\})|$, therefore $N(q_n)=V(P)\setminus(B(Q)\cup\{q_n\})$, that is, $q_n$ is adjacent to every
		vertex of $V(P)\setminus(B(Q)\cup\{q_n\})$.
	\end{proof} 
	
	A directed path $Q=q_1\cdots q_{\alpha}$ in $ D$ is \textbf{closed} if $(q_{\alpha},q_1)\in A(D)$. Equivalently, $Q$ is
	closed if and only if $q_1\cdots q_{\alpha}q_1$ is a Hamiltonian cycle of $D[Q]$.

\begin{figure}[htbp]
		\centering
		\begin{tikzpicture}[x=1cm,y=1cm,line cap=round,line join=round,
			lab/.style={font=\footnotesize},
			vlab/.style={font=\footnotesize,text depth=4pt},
			arr/.style={postaction={decorate,decoration={markings,
						mark=at position \spinepos\space with {\arrow[scale=1.0]{stealth}}}}},
			redge/.style={red,line width=0.85pt,shorten >=2.4pt,
				postaction={decorate,decoration={markings,
						mark=at position .5 with {\arrow[scale=1.0]{stealth}}}}},
			redln/.style={red,line width=1.1pt,shorten >=2.2pt}]
			\newcommand{\redln}[3]{%
				\draw[redln,postaction={decorate,decoration={markings,
						mark=at position #3 with {\arrow[scale=1.0]{stealth}}}}]
				(#1,0) -- (#2,0);}
			\def\pA{2.10}
			\def\pB{3.50}
			\def\pC{4.90}
			\def\pD{6.30}
			\def\pE{7.70}
			\def\pF{9.10}
			\def\pG{10.50}
			\def\pH{11.90}
			\def\pI{13.30}
			\def\pJ{14.70}
			\def\dv{0.12}
			\pgfmathsetmacro{\spinepos}{(14.00-2.10)/(14.70-2.10)}
			\draw[line width=0.8pt] (\pA,0) -- (\pJ,0);
			\redln{\pA}{\pD}{0.5}
			\redln{\pE}{\pI}{0.375}
			\draw[redge] (\pD,0) .. controls (8.568,1.200) and (12.432,1.200) .. (\pJ,0);
			\draw[redge] (\pJ,0) .. controls (12.367,0.773) and (10.033,0.773) .. (\pE,0);
			\foreach \x in {\pA,\pB,\pD,\pE,\pG,\pI,\pJ}{%
				\fill (\x,0) circle (1.2pt);}
			\foreach \x in {\pC,\pF,\pH}{%
				\node[lab,fill=white,inner sep=1.1pt] at (\x,0) {$\cdots$};}
			\node[vlab,anchor=north] at (\pA,-\dv) {$x_1$};
			\node[vlab,anchor=north] at (\pB,-\dv) {$x_2$};
			\node[vlab,anchor=north] at (\pD,-\dv) {$x_{j_P-1}$};
			\node[vlab,anchor=north] at (\pE,-\dv) {$x_{j_P}$};
			\node[vlab,anchor=north] at (\pG,-\dv) {$x_r$};
			\node[vlab,anchor=north] at (\pI,-\dv) {$x_{n-1}$};
			\node[vlab,anchor=north] at (\pJ,-\dv) {$x_n$};
			\node[lab,anchor=east] at (\pA-0.34,0) {$P'$};
		\end{tikzpicture}
		\caption{The path $P'$.}
		\label{fig:Pprime}
	\end{figure}
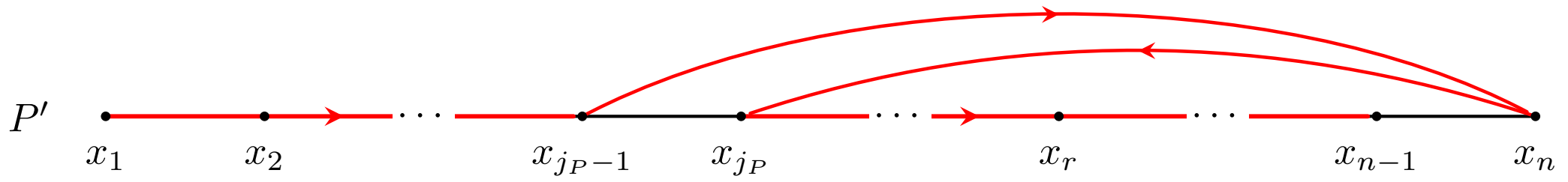

	\begin{lemma}\label{lem:shift}
		Let $s,r,K$ be positive integers with $s\geq r \geq 2$ and $s+r=K\geq 3$. Suppose that $D$ is an oriented graph with $\delta(D)\ge K$ and
		$d_D^+(u)+d_D^-(v)\ge K-1$ for every $(u,v)\notin A(D)$, but $D$ contains no $P(s,r)$. Suppose $P=x_1x_2\cdots x_n$ is a non-closed longest directed path of $D$, i.e., $(x_n,x_1)\notin A(D)$. 
		Then $P'$ is also a non-closed longest path with vertex set $V(P)$ (i.e., 
		$(x_{n-1},x_1)\not\in A(D)$), where
		\begin{equation}
			P':=x_1\cdots x_{j_P-1}\,x_n\,x_{j_P}\cdots x_{n-1}, \text{ where } j_P=\min\{i\in\{2,\dots,r\}:(x_n,x_i)\in A(D)\}.
			\label{eq:F}
		\end{equation}
	\end{lemma}
	\begin{proof}
		Assume that $r\ge2$ and that $P=x_1\cdots x_n$ is a longest path of $D$ with $(x_n,x_1)\notin A(D)$. By Lemma~\ref{longest path} (4), $(x_1,x_n)\in A(D)$ and
		$x_n$ is adjacent to $x_r$. If $(x_r,x_n)\in A(D)$, then the two paths
		$x_1\to\cdots\to x_r\to x_n$
		and
		$x_{n-s}\to\cdots\to x_n$
		have lengths $r$ and $s$, both end at $x_n$, and are internally disjoint as $r<n-s$. Hence
		$D$ contains $P(r,s)$, and so $P(s,r)$, a contradiction. Therefore $(x_n,x_r)\in A(D)$.
		
		Recall that  $j_P=\min\{i\in\{2,\dots,r\}:(x_n,x_i)\in A(D)\}.$
		If $j_P=2$, we have $(x_n,x_1)\not\in A(D)$. Now assume $j_P\geq 3$. 
		For every $2\le i<j_P$ we have $(x_n,x_i)\not\in A(D)$ by the minimality of $j_P$, by Lemma~\ref{longest path}(4), $(x_i,x_n)\in A(D)$ for $2\le i<j_P$.
		By the selection of $j_P$, we have $P':=x_1\cdots x_{j_P-1}\,x_n\,x_{j_P}\cdots x_{n-1}$ is a path with vertex set $V(P)$ (see Figure~\ref{fig:Pprime}), so $P'$ is a longest path. 
		
		The following proof is divided into two cases to show $(x_{n-1},x_1)\not\in A(D)$.
		
		\smallskip
		\noindent\textit{Case 1: $n-s-r-1=0$.} 
		Then $R(P')=B(P')=\emptyset$. By
		Lemma~\ref{longest path}(4), the terminal vertex of every longest path on $V(P)$ is adjacent to all other vertices of $V(P)$.
		If $D[V(P)]$ contains a Hamilton cycle $C$, then $C$ can be cut at any vertex, so every vertex of $V(P)$ is the terminal vertex of a longest path. By Lemma~\ref{longest path}(4), the terminal vertex is adjacent to all the other vertices. Thus $D[V(P)]$ is a tournament.
		Since $s\ge r\ge2$, the path $P(s,r)$ is not antidirected, so by Theorem~\ref{Havet}, $D[V(P)]$ contains $P(s,r)$, a contradiction. 
		Hence $D[V(P)]$ contains no Hamilton cycle, and thus in $P'$, $(x_{n-1},x_1)\not\in A(D)$ as desired.
		
		\smallskip
		\noindent\textit{Case 2: $n-s-r-1>0.$}
		Since $P=x_1\cdots x_n$ is a longest path of $D$ with $(x_n,x_1)\notin A(D)$, by Lemma~\ref{longest path} (4), $(x_1,x_n)\in A(D)$ and
		$x_n$ is adjacent to $x_r$.  Now Lemma~\ref{longest path}(3) gives $s=r$. If $(x_{n-1},x_1)\in A(D)$, then $x_1\cdots x_{{j_P}-1}\,x_n\,x_{j_P}\cdots x_{n-1}x_1$ is a cycle. Deleting the arc $x_nx_{j_P}$ from the cycle gives the longest path
		$$
		P''=x_{j_P}\cdots x_{n-1}x_1\cdots x_{{j_P}-1}x_n,
		$$
		in which $x_{r+1}\notin B(P'')$, so $x_n$ and $x_{r+1}$ are adjacent by Lemma~\ref{longest path}(4). But in $P$ the vertex $x_{r+1}$ lies in $B(P)$, so they are nonadjacent by Lemma~\ref{longest path}(1), a contradiction. Hence $(x_{n-1},x_1)\not\in A(D)$.
		
		Both cases give a contradiction, which completes the proof of the Lemma.
	\end{proof}
	
	\subsection{Completion of the Proof}
\begin{figure}[htbp]
		\centering
		\begin{tikzpicture}[x=1cm,y=1cm,line cap=round,line join=round,
			lab/.style={font=\small},
			vlab/.style={font=\small,text depth=5pt},
			elab/.style={font=\small,text depth=5pt},
			arr/.style={postaction={decorate,decoration={markings,
						mark=at position 1 with {\arrow[scale=1.25]{stealth}}}}}
			]

			\def\yA{0}      \def\yB{-1.20}   \def\yC{-2.40}
			
			\def\pA{2.10}    
			\def\pB{3.45}    
			\def\pC{4.80}    
			\def\pD{6.15}    
			\def\pE{7.55}    
			\def\pK{9.69}    
			\def\pF{10.35}   
			\def\pG{11.75}   
			\def\pH{12.75}   
			\def\pI{14.35}   
			
			\def\xa{\pA}     
			\def\xb{\pI}     
			\def\dv{0.14}

			\newcommand{\prow}[2]{%
				\draw[line width=1.0pt,arr] (\xa,#1) -- (\xb,#1);
				\foreach \x in {\pA,\pB,\pD,\pE,\pF,\pH,\pI}{%
					\fill (\x,#1) circle (1.5pt);}
				\foreach \x in {\pC,\pK,\pG}{%
					\node[lab,fill=white,inner sep=1.2pt] at (\x,#1) {$\cdots$};}
				\node[lab,anchor=east] at (1.34,#1) {#2};}
			
			\prow{\yA}{$P^{(0)}$}
			\prow{\yB}{$P^{(1)}$}
			\prow{\yC}{$P^{(2)}$}

			\node[vlab,anchor=north] at (\pA,\yA-\dv) {$x_1^{(0)}$};
			\node[vlab,anchor=north] at (\pB,\yA-\dv) {$x_2^{(0)}$};
			\node[vlab,anchor=north] at (\pD,\yA-\dv) {$x_h^{(0)}$};
			\node[vlab,anchor=north] at (\pE,\yA-\dv) {$x_j^{(0)}$};
			\node[vlab,anchor=north] at (\pF,\yA-\dv) {$x_r^{(0)}$};
			\node[vlab,anchor=north] at (\pH,\yA-\dv) {$x_{n-1}^{(0)}$};
			\node[vlab,anchor=north] at (\pI,\yA-\dv) {$x_n^{(0)}$};

			\node[vlab,anchor=north] at (\pA,\yB-\dv) {$x_1^{(1)}$};
			\node[vlab,anchor=north] at (\pB,\yB-\dv) {$x_2^{(1)}$};
			\node[vlab,anchor=north] at (\pD,\yB-\dv) {$x_h^{(1)}$};
			\node[elab,anchor=north west] at (\pE-0.48,\yB-\dv) {$x_j^{(1)}=x_n^{(0)}$};
			\node[vlab,anchor=north] at (\pF,\yB-\dv) {$x_r^{(1)}$};
			\node[vlab,anchor=north] at (\pH,\yB-\dv) {$x_{n-1}^{(1)}$};
			\node[elab,anchor=north] at (\pI,\yB-\dv) {$x_n^{(1)}=x_{n-1}^{(0)}$};

			\node[vlab,anchor=north] at (\pA,\yC-\dv) {$x_1^{(2)}$};
			\node[vlab,anchor=north] at (\pB,\yC-\dv) {$x_2^{(2)}$};
			\node[vlab,anchor=north] at (\pD,\yC-\dv) {$x_h^{(2)}$};
			\node[elab,anchor=north west] at (\pE-0.48,\yC-\dv) {$x_{j'}^{(2)}=x_n^{(1)}$};
			\node[vlab,anchor=north] at (\pF,\yC-\dv) {$x_r^{(2)}$};
			\node[vlab,anchor=north] at (\pH,\yC-\dv) {$x_{n-1}^{(2)}$};
			\node[elab,anchor=north] at (\pI,\yC-\dv) {$x_n^{(2)}=x_{n-1}^{(1)}$};

			\pgfmathsetmacro{\xmid}{(\xa+\xb)/2 - 1.43}
			\node[lab,anchor=north] at (\xmid,\yC-0.67) {$\vdots$};
			
		\end{tikzpicture}
		\caption{Two steps of the iteration.}
		\label{fig:iteration3}
	\end{figure}
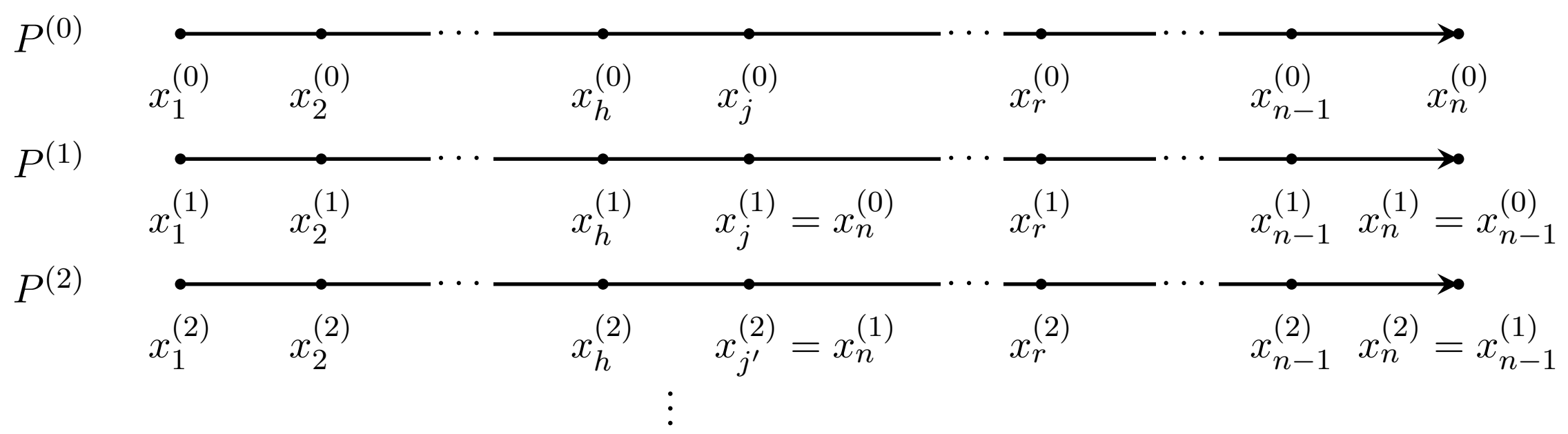

	\begin{proof}[Proof of Theorem~\ref{thm:main}]

		Suppose, for a contradiction, that an oriented graph $D$ satisfies $\delta(D)\ge s+r$ and $d_D^+(u)+d_D^-(v)\ge s+r-1$ for every $(u,v)\notin A(D)$, but contains no $P(s,r)$. At the beginning of this section, we assumed that $s\ge r\ge1$, $s\ge 2$ and $s+r\ge3$. Let $P=x_1\cdots x_n$ be a longest path of $D$.

		The proof proceeds by induction on $r$. 
		Let $r=1$.
		By Lemma~\ref{longest path}(3) we have
		$R(P)=B(P)=\emptyset$, by (4), $x_n$ is adjacent to every other vertex of $V(P)$. If $(x_1,x_n)\in A(D)$,  $x_2\to\cdots\to x_n \leftarrow x_1$, is $P(s,1)$. Otherwise $(x_n,x_1)\in A(D)$, will be a directed Hamilton cycle on $V(P)$. Cutting this cycle at any vertex and using Lemma~\ref{longest path}(4) again, $D[V(P)]$ is a tournament containing the non-antidirected path $P(s,1)$ by Theorem \ref{Havet}.
		Assume that for every oriented graph $D$ with $\delta(D)\ge s+f$ and $d_{D}^+(u)+d_{D}^-(v)\ge s+f-1$,  $D$ contains $P(s,f)$ for all $f<r$.

		Now we prove that $D$ contains $P(s,r)$, where $r\ge2$, by contradiction. 
		So assume that $D$ contains no $P(s,r)$. Indeed, $P$ is a  non-closed longest path of $D$. For otherwise, $D[V(P)]$ contains a Hamilton cycle $C$, and then $C$ can be cut at any vertex, so every vertex of $V(P)$ is the terminal vertex of a longest path. By Lemma~\ref{longest path}(4), the terminal vertex is adjacent to all the other vertices. Thus $D[V(P)]$ is a tournament.
		Since $s\ge r\ge2$, the path $P(s,r)$ is not antidirected, so by Theorem~\ref{Havet}, $D[V(P)]$ contains $P(s,r)$, a contradiction.
		
		For each non-closed longest path, assume $Q:=q_1\cdots q_n$, of $D$, Lemma~\ref{lem:shift} always  provides a new non-closed longest path $Q':=q_1\cdots q_{j_Q-1}\,q_n\,q_{j_Q}\cdots q_{n-1}$ on $V(Q)$ that is also non-closed, where $j_Q=\min\{i\in\{2,\dots,r\}:(q_n,q_i)\in A(D)\}$. Consequently, applying the replacement of the non-closed longest path of $D$ as stated in Lemma \ref{lem:shift}  repeatedly, starting from $P$, we obtain a sequence of non-closed longest paths
		$$
		P:=P^{(0)}\ \longrightarrow\ P^{(1)}\ \longrightarrow\ P^{(2)}\ \longrightarrow\ \cdots .
		$$
		The iteration terminates at the first time that a newly produced non-closed longest path $P^{(t)}$ appears in the collection $\{P^{(i)}:i\in [t-1]\}$ (see Figure~\ref{fig:iteration3}). This iterative procedure always terminates, since $n$ is finite: there are only finitely many distinct paths on an $n$-vertex set. Hence, assume the iterative procedure terminates at step $t$, and output the sequence of non-closed longest paths generated successively throughout the iteration $\mathcal{P}:=\{P^{(i)}:=x^{(i)}_1 \  x^{(i)}_2\ \cdots \,x^{(i)}_n:i\in [t-1]\}$. Observe that the replacement operation for non-closed longest paths given by Lemma~\ref{lem:shift} yields
		\begin{equation}\label{pi}
			P^{(i)}:=x^{(i-1)}_1\cdots x^{(i-1)}_{j_{P^{(i-1)}}-1}\,x^{(i-1)}_n\,x^{(i-1)}_{j_{P^{(i-1)}}}\cdots x^{(i-1)}_{n-1},
		\end{equation}
		where $j_{P^{(i-1)}}=\min\{\ell \in\{2,\dots,r\}:(x^{(i-1)}_n,x^{(i-1)}_{\ell})\in A(D)\}$.

		Let $W=V(P)\setminus  \{x^{(i)}: i\in [0,t-1]\} $ be the set of vertices that are never the terminal vertex during the process. 
		Indeed, $W\neq \emptyset$, it follows from $j_{P^{(i)}}\geq 2$ for all $i\in [0,t-1]$ that $x_1\in W$.
		Assume $h=|W|$, so $h<n$. We assert that $W:=\{x_1,x_2,\ldots,x_h\}$, i.e., the set of the first $h$ vertices on the path $P$.  For otherwise, there is a vertex $x_i\notin W$ with $x_{i+1}\in W$. The definition of $W$ implies that  $x_i$ is the terminal vertex of some path $P^{(j)}$, where $j\in [0,t-1]$. And \eqref{pi} implies that $P^{(j-1)}$ terminates at $x_{i+1}$, i.e., $ x_{i+1}\in \{x^{(i)}: i\in [0,t-1]\} $ , which contradicts $x_{i+1}\in W$.

		Put $U=\{x^{(i)}: i\in [0,t-1]\}$. 
		Since $x_n\in U$,  $U\ne\emptyset$.
		
		\begin{claim}
			We obtain that $W\to U$.
		\end{claim}
		\begin{proof}
			Let $u\in U$. 
			By the definition of $U$, the vertex $u$ is the terminal vertex of some longest path $P^{(a)}$ on $V(P)$ occurring in the sequence $\mathcal{P}$. 
			And the insertion position $j_{P^{(a)}}=\min\{\ell \in\{2,\dots,r\}:(x^{(a)}_n,x^{(a)}_{\ell})\in A(D)\}$ satisfies $j_{P^{(a)}}>h$: $W=V(P)\setminus  \{x^{(i)}: i\in [0,t-1]\} =\{x_1,x_2,\ldots,x_h\}$ is the first $h$ vertices of $P$. Hence $ \{x_1,x_2,\ldots,x_h\}$ are also the first $h$ vertices of $P^{(a)}$. So $x_i\notin B(P^{(a)})$ for all $i\in [h]$. 
			By the definition of $j_{P^{(a)}}$, the vertex $u$ dominates none of the first $j_{P^{(a)}}-1$ vertices, and hence $u\not\to x_i$ for every $i\in [h]$. This forces
			$x_i\to u$. Therefore $W\to U.$     
		\end{proof}

		\begin{claim}
			The oriented graph $D[U]$ satisfies $\delta(D[U])\ge s+r-h$ and $d^+_{D[U]}(u)+d^-_{D[U]}(v)\ge s+r-h-1$ for every $(u,v)\notin A(D[U])$.
		\end{claim}
		\begin{proof}
			By Lemma~\ref{longest path}(2) applied to the longest path terminating at $u$, and by the maximality of that path, the vertex $u$ has no neighbour outside $V(P)$. 
			That is, $N_D(u)\subseteq V(P)$. By Lemma~\ref{longest path}(3), $d_D(u)=s+r$.
			Since $V(P)=W\cup U$ and every vertex of $W$ sends an arc to $u$, by $W\to U$ we obtain
			$$
			d_{D[U]}(u)=d_D(u)-h=s+r-h.
			$$
			As $u\in U$ was arbitrary, this yields
			\begin{equation}
				\delta(D[U])\ge s+r-h. \label{eq:delta-U}
			\end{equation}
			
			Now let $u,v\in U$ be distinct with $(u,v)\notin A(D)$. 
			By $W\to U$ every vertex of $W$ sends an arc both to $u$ and to $v$, whereas $u$ sends no arc to $W$. 
			Hence $d_D^+(u)=d_{D[U]}^+(u)$, while the in-neighbours of $v$ in $D$ are those in $U$ together with the $h$ vertices of $W$, so $d_D^-(v)=d_{D[U]}^-(v)+h$. Substituting these two identities into
			the degree condition $d_D^+(u)+d_D^-(v)\ge s+r-1$ gives
			\begin{equation}
				d_{D[U]}^+(u)+d_{D[U]}^-(v)\ge s+r-h-1. \label{eq:missing-U}
			\end{equation}
			By \eqref{eq:delta-U} and \eqref{eq:missing-U}, the oriented graph $D[U]$ satisfies
			$\delta(D[U])\ge s+r-h$ and
			$d_{D[U]}^+(u)+d_{D[U]}^-(v)\ge s+r-h-1$ for every $(u,v)\notin A(D[U])$, as desired.
		\end{proof}

		By the induction hypothesis, $D[U]$ contains $P(s,r-h)$. 
		Let $a_0\to a_1\to\cdots\to a_{r-h}$ be its path of length $r-h$, the second path, which has length $s$, also ends at $a_{r-h}$.
		Since $a_0\in U$, $W\to U$ gives
		$x_h\to a_0$. 
		Hence
		$$
		x_1\to\cdots\to x_h\to a_0\to\cdots\to a_{r-h}
		$$
		is a directed path of length $(h-1)+1+(r-h)=r$. 
		It is internally disjoint from the path of length $s$ as $W\cap U=\emptyset$. 
		So $D$ contains $P(r,s)$, i.e.,  $P(s,r)$, a contradiction.
		
		It remains to obtain $Q(s,r)$. Let $D^{\mathrm{rev}}$ be the reverse of $D$.
		Reversing every arc preserves the degree condition: $\delta(D^{\mathrm{rev}})
		=\delta(D)$, and if $(u,v)\notin A(D^{\mathrm{rev}})$ then $(v,u)\notin A(D)$, so
		$d_{D^{\mathrm{rev}}}^+(u)+d_{D^{\mathrm{rev}}}^-(v)=d_D^-(u)+d_D^+(v)\ge s+r-1$.
		Hence $D^{\mathrm{rev}}$ contains $P(s,r)$, and reversing this back yields $Q(s,r)$ in $D$. This completes the proof.
	\end{proof}

	\section{Conclusion}
	
	Jackson~\cite{Jackson1981} proved the one-block case (directed paths) under the
	Ore-type condition.
	In this paper we settle the two-block case under the same condition. 
	A natural open problem is whether the analogue holds for oriented paths with three or more blocks.

	\section*{Declaration of generative AI}
	During the development and preparation of this work, the authors used ChatGPT as an auxiliary tool to assist in checking preliminary mathematical arguments, identifying potential gaps in proofs, and exploring possible approaches to addressing them. The authors critically evaluated the suggestions and substantially revised the relevant arguments through multiple rounds of mathematical examination and rewriting. The authors take full responsibility for the correctness, originality, and final content of the manuscript.

\end{document}